\documentclass[11pt]{amsart}
\makeatletter

\title[Hilbert series of fiber cones of ideals]{ Hilbert Series of fiber 
cones of ideals with almost minimal
mixed multiplicity}
\thanks{ {\it Key words }: Hilbert series, Cohen-Macaulay rings, almost
minimal mixed multiplicity, mixed
multiplicities, fiber cones, reduction of an ideal, joint reduction.\\
 Presented in the second 
{\it National Meeting in Commutative Algebra and Algebraic Geometry }
held in February 1999 at Harish-Chandra Research Institute, Allahabad,
India.
}

\author{  Clare D'Cruz}
\address{Chennai Mathematical Institute, G. N. Chetty Road, T. Nagar,
Chennai 600 017 India }
\email {clare@cmi.ac.in}
\author{J. K. Verma}
\address{Department of Mathematics, Indian Institute of Technology, Bombay,
Mumbai 400 076 India} 
\email{ jkv@math.iitb.ac.in }

\newcommand{\ncom}{\newcommand}
\ncom{\beqn}{\begin{eqnarray*}}
\ncom{\eeqn}{\end{eqnarray*}}
\ncom{\beq}{\begin{eqnarray}}
\ncom{\eeq}{\end{eqnarray}}
\ncom{\been}{\begin{enumerate}}
\ncom{\eeen}{\end{enumerate}}
\ncom{\nno}{\nonumber}
\ncom{\hs}{\mbox{\hspace{.25cm}}}
\ncom{\rar}{\rightarrow}
\ncom{\lrar}{\longrightarrow}
\ncom{\Rar}{\Rightarrow}
\ncom{\noin}{\noindent}
\newtheorem{thm}{Theorem}[section]
\newtheorem{lemma}[thm]{Lemma}
\newtheorem{cor}[thm]{Corollary}
\newtheorem{pro}[thm]{Proposition}
\newtheorem{example}[thm]{Example}
\newtheorem{definition}[thm]{Definition}
\newtheorem{remark}[thm]{Remark}
\ncom{\bt}{\begin{thm}}
\ncom{\et}{\end{thm}}
\ncom{\bl}{\begin{lemma}}
\ncom{\el}{\end{lemma}}
\ncom{\bco}{\begin{cor}}
\ncom{\eco}{\end{cor}}
\ncom{\bp}{\begin{pro}}
\ncom{\ep}{\end{pro}}
\ncom{\bex}{\begin{example}}
\ncom{\eex}{\end{example}}
\ncom{\bd}{\begin{definition}}
\ncom{\ed}{\end{definition}}
\ncom{\brm}{\begin{remark}}
\ncom{\erm}{\end{remark}}
\ncom{\comx}{I\!\!\!\!C}
\ncom{\proj}{I\!\!\!P}
\ncom{\zee}{$Z\!\!\!\!Z$}
\ncom{\ze}{Z\!\!\!\!Z}
\ncom{\Q}{$I\!\!\!\!Q$}
\ncom{\N}{I\!\!N}
\ncom{\sz}{\scriptsize}
\ncom{\CM}{Cohen-Macaulay }
\ncom{\sop}{system of parameters}
\ncom{\eop}{\hfill{$\Box$}}
\ncom{\tfae}{the following are equivalent:}
\ncom{\mm}{minimal multiplicity }
\ncom{\f}{\frac}
\ncom{\la}{\lambda}
\ncom{\si}{\sigma}
\ncom{\ssize}{\scriptsize}
\ncom{\al}{\alpha}
\ncom{\be}{\beta}
\ncom{\Si}{\Sigma}
\ncom{\ga}{\gamma}
\ncom{\kbar}{\overline{\kappa}}
\ncom{\bib}{\bibitem}
\ncom{\sst}{\subset}
\ncom{\sms}{\setminus}
\ncom{\seq}{\subseteq}
\ncom{\est}{\emptyset}
\ncom{\pf}{{\bf Proof: }}
\ncom{\bighs}{\hspace{.5 cm}}
\ncom{\ulin}{\underline}
\ncom{\olin}{\overline}
\ncom{\bip}{\bigoplus}
\ncom{\sta}{\stackrel}
\begin{document} 
\maketitle

\begin{center}
{\it Dedicated to Bill Heinzer on the occasion of his sixtieth birthday}
\end{center}

\begin{abstract} Let $(R, m)$ be a \CM local ring and $I$ be an
$m$-primary ideal. We introduce ideals of almost minimal mixed
multiplicty which are  analogues of ideals studied by J. Sally in
\cite{s}. The main theorem describes the Hilbert series of fiber cones of
these ideals.

\end{abstract} 

\section{Introduction}
Throughout this paper $R$ will denote a $d$-dimensional Noetherian local
ring with
unique maximal ideal $m$ and $I$ will denote an $m$-primary ideal. The
fiber cone of $I$, $F(I)$, is defined to be the graded ring
$\bigoplus_{n \geq 0}I^n/mI^n.$ The fiber cone is the homogeneous
coordinate ring of the fiber over the closed point $m$ in the blowup
of
Spec\,$R$ along the subscheme  Spec\,$(R/I).$ The fiber cone is in the
class of rings called the blowup rings associated with $I$. Recently it
has been investigated by a number of researchers (  \cite{cz1},
 \cite{cz2}, \cite{dgh}, \cite{drv}, \cite{g}, \cite{hh} \cite{sh1},
\cite{sh2}).

In the last few decades the Rees algebra and the associated graded ring of
an ideal (defined below) have been investigated by many researchers.
However the fiber cone of an ideal has not been studied much. The fiber
cones are useful in a number of situations. For example B. Singh
characterized prime ideals $p$ that can be permissible centers of blowing
up in terms of the Hilbert series of the fiber cone of $p$ \cite{si}. 
Fiber cones turn out to be useful in understanding evolutions( \cite{em},
\cite{hh}). In the latter case the Cohen-Macaulay property of the fiber
cone is useful while in the former case one needs to have prior knowledge
of the Hilbert series of the fiber cone.  

Our objective in this paper is to calculate Hilbert series of fiber cones
of ideals 
of almost minimal mixed multiplicity introduced below. These ideals
are analogues for  $m$-primary ideals of ideals studied by J. Sally in
\cite{s}. The knowledge of Hilbert series of $F(I)$ provides us useful
information about the number of generators of all powers of $I$ and more
importantly  
in view of the main theorem of \cite{drv} this  knowledge   helps in
detecting its Cohen-Macaulay property. 

Let $\mu(I)$ denote the minimum number of
generators of an ideal $I.$  The Hilbert series of $F(I)$ is defined by 
$H(F(I),\la):=\sum_{n=0}^{\infty}\mu(I^n)\la^n .$ 
 In order to state the main theorem of this paper we recall the necessary
definitions
first. Let $\ell(.)$ denote length. It was proved in \cite{b} that  for
large values of $r$ and $s,$ the
function $\ell(R/m^rI^s)$ is given by a polynomial $P(r,s)$ of total
degree $d$ in $r$ and $s$. The polynomial
$P(r,s)$ can be written in the form:
$$ P(r,s)=\sum_{i+j \leq d}e_{ij}\binom{r+i} {i}\binom{s+j}{ j},$$
\noindent
where $e_{ij}$ are certain integers. When $i+j=d,$ we put
$e_{ij}=e_j(m|I)$
for $j=0,1, \ldots,d.$ In this case these are called the mixed
multiplicities of $m$ and $I.$ It is known that $e_0(m|I)= e(m)$ and
$e_d(m|I)= e(I)$ where $e(.)$ denotes the Hilbert-Samuel
multiplicity \cite{r1} .  The other mixed multiplicities too can be shown
to be
Hilbert-Samuel multiplicities of certain systems of parameters (\cite{t}
and \cite{r2}).  We recall an important result from \cite{r2}. In this
important paper  Rees introduced joint reductions to calculate mixed 
multiplicities. A
set of elements $x_1, \ldots, x_d$ is a called a {\em joint reduction} of
a set of ideals $I_1, \ldots, I_d$ if $x_i \in I_i$ for $i=1, \ldots, d$
and there exists a positive integer $n$ so that \beqn
  \left[ \sum_{j=1}^{d} x_j \, I_1 \cdots \hat{I_j} \cdots I_d \right]
   (I_1 \cdots I_d)^{n-1}
=  (I_1 \cdots I_d)^n.
\eeqn
Rees proved that if $R/m$ is infinite, then joint reductions
exist and $e_j(m|I)$ is the multiplicity of any joint reduction of  
the multiset of ideals consisting of $j$ copies of $I$ and $d-j$ copies of
$m.$
We shall denote such a multiset  of ideals by $(m^{[d-j]}|I^{[j]})$.

It has been proved in \cite{drv} that if $R$ is Cohen-Macaulay then
$e_{d-1}(m|I) \geq \mu(I)-d+1.$  We say that $I$ has minimal mixed
multiplicity if
$e_{d-1}(m|I)= \mu(I)-d+1.$ We calculated the Hilbert series
of $F(I)$ when $I$ has minimum mixed multiplicity and
showed in \cite{drv} that for these ideals $F(I)$ is \CM if and only if
the reduction
number of $I$ is at most one. See the next section  for the
definition of reduction number of an ideal. 

\bd  
Let $(R,m)$ be a \CM local ring of dimension $d$. An $m$-primary
ideal $I$ of $R$ is said to have almost minimal mixed multiplicity if
$e_{d-1}(m|I)= \mu(I)-d+2.$ 
\ed

\bd
Let $J=(x, a_1, a_2, \ldots, a_{d-1}) $ be a joint reduction
of $(m|I^{[d-1]}).$ Define $r_J(m|I)$ to be the smallest integer $n$, if
it exists, so that $mI^n=xI^n+(a_1, a_2, \ldots, a_{d-1})mI^{n-1}$. 
The smallest of all $r_J(m|I)$ where $J$ is varying is denoted by
$r(m|I)$ If
there is no such
integer we say that $r_J(m|I)$ is infinite and we write $r(m|I)=\infty.$
\ed  

 Let $\gamma(I)$ and
$\phi(I)$ denote the depths of the ideals generated by elements of
positive degree in $G(I):= \bigoplus_{n \geq 0}I^n/I^{n+1}$
and
$F(I)$ respectively. We can now state the
main theorem of this paper:

\bt
Let $(R,m)$ be a \CM local ring of dimension $d$ with infinite residue
field. Let $I$ be an $m$-primary ideal of almost minimal mixed
multiplicity. Let  $\gamma(I) \geq d-1$ and  $\phi(I) \geq d-2.$ Then  the
Hilbert series of $F(I)$ is given by

\beqn
H(F(I), \la)
= \left\{     
         \begin{array}{ll}
         \displaystyle \f{1+ (\mu(I)-d)\la}{(1 - \la)^d}
          &  \,\, \mbox{if} \,\, r(m|I)=\infty \\ \\
          \displaystyle \f{1 + (\mu(I)-d ) \la  + \la^{r(m|I)}}{(1 -
           \la)^d}
          &  \,\, \mbox{if} \,\, r(m|I) \,\,\mbox{is finite.}  \\ \\
         \end{array}\right.
\eeqn
\et

By the results proved in \cite{drv}, the following corollary is easily
deduced.
\bco
Under the conditions of the above theorem, $F(I)$ is \CM if and only if
either $r(I) \leq 1$ or $r(I)=2$ and $\ell(I^2/JI+mI^2)=1$ for some(and
hence all) 
minimal reduction $J$ of $I.$
\eco

The paper is organized as follows: In Section 2 we characterize ideals
of minimal mixed multiplicity. The main  theorem requires a different
approach in dimension one and therefore we have proved it in Section 3.
In Section 4, we apply basic results about mixed multiplicities of ideals
to find the Hilbert series of fiber cones of ideals of almost minimal
mixed multiplicity in the two dimensional case. In Section 5 the main
theorem is proved for all
Cohen-Macaulay local rings of dimension $\geq 2.$ In section 6 we provide
some examples to illustrate the main theorem. Finally we answer a question
about Cohen-Macaulay property of fiber cones raised in \cite{dgh}.

\noindent
{\bf Acknowledgement}: The first author thanks the Indian Institute of
Technology, Bombay for its hospitality where this work was completed. 
The second author thanks Dan Katz for useful
discussions while visiting the University of Kansas. The authors also thank the
referee for a careful reading and  pointing out the references \cite{hh} and \cite{em}.


\section{Preliminaries}
\noindent
In this section we  prove some preliminary results which will be used in
the
subsequent sections. We begin by characterizing ideals of almost
minimal mixed multiplicity.

\bp
\label{char}
Let $(R,m)$ be a \CM local ring with infinite residue field.
Suppose $I$ is an $m$-primary ideal. Then $I$ has almost minimal mixed
multiplicity if and only if for any joint reduction 
$(x, a_1, a_2, \ldots, a_{d-1})$ of $(m|I^{[d-1]})$, 
$$\ell(mI/xI+  (a_1, a_2, \ldots, a_{d-1})m)=1.$$
\ep

\pf From the proof of Theorem 2.4 of \cite{drv} we
have 
\beqn
   \ell \left( \f{R}{xI + (a_{1}, \ldots, a_{d-1}) m} \right)  
 - \ell \left( \f{R}{I} \right)
=   d-1 + e_{d-1}(m|I).
\eeqn

Hence
\beqn  \ell \left( \f{mI}{xI + (a_{1}, \ldots, a_{d-1}) m} \right)  
       +\mu(I)
= d-1+ e_{d-1}(m|I).
\eeqn

The proposition is clear from the above equation.

\bl
\label{length}
Let $(R,m)$ be a \CM local ring of positive dimension $d$. Let $I$ be an 
$m$-primary ideal of $R$ having almost minimal mixed multiplicity. Then for
any joint reduction $(x, a_1, \ldots, a_{d-1})$  of $(m|I^{[d-1]}),$

$$
\ell(mI^n / (a_1, \ldots, a_{d-1}) mI^{n-1} +  xI^n) \leq 1  
\,\,\, \mbox{for all}\,\,\, n \geq 1.$$

\el

\pf  Put $J = (a_1, \ldots, a_{d-1})$.
Since $\mu(I) = e_{d-1}(m|I) + d-2$, by Proposition \ref{char} $\ell
(Im/xI +Jm)= 1.$
Hence  there exists $y \in m$ and $b \in I$ such 
that $mI = Jm + xI + (yb)$ and $myb \subset Jm + xI .$ We claim that for
all $n \geq 1,$ $ mI^n = JmI^{n-1} + xI^n + (yb^n).$ 
If $n = 1$, we are done by Proposition \ref{char}. If $n >1$, then by
induction hypothesis
$$  
mI^n  
= (mI^{n-1}) I 
= (JmI^{n-2} + xI^{n-1} + (y b^{n-1})) I
= JmI^{n-1} + xI^n+ (y b^{n-1})I. 
$$
\noindent
Since $$
(y b^{n-1}) I \seq
             bmI^{n-1}
           = (JmI^{n-2}  + xI^{n-1}+ (y b^{n-1}))b
         \seq JmI^{n-1} + xI^n + (y b^{n}),
$$
$mI^n  = JmI^{n-1} + xI^n + (y b^n).$
Since $m (yb) \seq Jm+xI$, 
$$
m(yb^n) = m(yb)(b^{n-1}) \seq (Jm + xI)(b^{n-1}) \seq JmI^{n-1} + xI^n. 
$$
Hence $\ell(mI^n/JmI^{n-1} + xI^n) \leq 1$.

\noindent
{\bf Notation: } Let $I$ be an ideal of a ring $R$. For an element $a \in
I$ let $a^*$ (resp. $a^o$) denote its residue class  in $I/mI$ 
(resp. $I/I^2$).

\bl
\label{joint}
Let $(R,m)$ be a \CM local ring of dimension $d$ and $I$ an $m$-primary
ideal of $R.$ Suppose $\gamma(I) \geq d-1$ and $\phi(I) \geq d-2.$
Then there exists a joint reduction $(x, a_1, a_2, \ldots, a_{d-1})$ of
$(m|I^{[d-1]})$ such that $a_1^o, \ldots, a_{d-1}^o$ and $a_1^*, \ldots,
a_{d-2}^*$ are regular sequences in $G(I)$ and $F(I)$ respectively.  
\el
\pf 
By Lemma 1.2 of \cite{r2} and its consequences, we can find the desired  
joint reduction by avoiding the contractions to $R$  of the associated
primes of $(t^{-1},m)$ and  $(t^{-1})$ in $R[It,t^{-1}]$ at each stage 
and using Lemma  \ref{fiberpreserve}.


\section{Main theorem in dimension one}
In this section we will prove the main theorem for one dimensional
\CM local rings. We first recall the concept of reduction of ideals from
\cite{nr}. Let $(R,m)$ be a local ring with infinite residue field. An
ideal $J$ of $R$ is called a reduction of an ideal $I$ of $R$ if $J\subset
I$ and there exists an integer $n$ such that $JI^n=I^{n+1}.$ If $J$ is
the smallest such ideal then it is called a minimal reduction of $I.$ All
minimal reductions of $I$ are generated by the same number of elements
which is
called the analytic spread of $I$ and it is equal to the Krull dimension
of $F(I).$ The reduction number $r_J(I)$ of an ideal $I$ with respect to a
minimal reduction $J$ is the smallest integer $n$ such that
$JI^n=I^{n+1}.$ The reduction number $r(I)$ of $I$ is the minimum of
$r_J(I)$ where $J$ ranges over all the minimal reductions of $I.$

 Throughout this section $(R,m)$
will denote a \CM local ring of dimension one with infinite residue field.

\bl
\label{hser1}
 Let $(x)$ be a minimal reduction of $m$. Then 
$$
H(F(I), \la) = \f{e(R)}{1 - \la}
             - H \left( \f{mR[It]}{xR[It]}, \la \right).
$$
\el
\pf Since $(x)$ is a minimal reduction  of $m$, $x$ is a
nonzerodivisor.  Therefore $\ell(xR/xI^n) = \ell(R/I^n)$ for all  $n
\geq 0$. Hence for all $n \geq 0,$

\beq
\label{one}
\mu(I^n) 
= \ell \left( \f{R}{xR} \right)
+ \ell \left( \f{xR}{xI^n} \right)
- \ell \left( \f{mI^n}{xI^n} \right)
- \ell \left( \f{R}{I^n} \right)
=  e(R) 
- \ell \left( \f{mI^n}{xI^n} \right).
\eeq
Therefore
\beqn
   H(F(I), \la)
= \sum_{n=0}^{\infty} \mu(I^n) \la^n
=  \sum_{n \geq 0} 
    \left[ e(R) - \ell \left( \f{mI^n}{xI^n} \right) \right] \la^n
= \f{e(R)}{1 - \la} 
- H \left( \f{mR[It]}{xR[It]}, \la \right).
\eeqn

\bp
 Suppose $\mu(I) = e(R)$. Then 
$$
H(F(I), \la)
= \f{1 + (e(R) -1) \la}{1 - \la}.
$$
\ep
\pf We have $\mu(I) = e(R)$ if and only if $mI=xI$ by (\ref{one}).
 Hence $mI^{n}
= xI^{n}$ for all $n \geq 1$. Therefore by Lemma~\ref{hser1},
\beqn
H(F(I), \la) 
= \f{e(R)}{1 - \la} - \ell \left( \f{m}{xR} \right)
= \f{e(R)}{1 - \la} - e(R) + 1
= \f{1 + (e(R) -1) \la}{1 - \la}.
\eeqn

\bt
\label{main1}
  Suppose $\mu(I) = e(R) - 1$. Then  the Hilbert series of $F(I)$ is given
by

\beqn
H(F(I), \la) 
= \left\{ 
         \begin{array}{ll}
         \displaystyle \f{1+ (e(R) - 2)\la}{1 - \la} 
          &  \mbox{if}\,\,\, r(m|I)=\infty  ,\\ \\
          \displaystyle \f{1 + (e(R) - 2) \la  + \la^s}{1 - \la}
          & \;\;  \mbox{if}\,\,\,   r(m|I)=s. \\ 
\end{array} 
\right.
\eeqn

\et

\pf Let $(x)$ be a minimal reduction of $m$. Apply Lemma~\ref{hser1}.
If $\ell (mI^n/xI^n) =1$ for all $n
\geq 1$ then 
\beqn
  H(F(I), \la)
&=& \f{e(R)}{1- \la} - \ell \left( \f{m}{xR}\right)
  - \sum_{n \geq 1} \ell \left( \f{mI^n}{xI^n}\right) \la^n\\
&=& \f{e(R)}{1 - \la} - (e(R)-1) - \sum_{n \geq 1} \la^n\\
&=& \f{1 + (e(R) - 1) \la }{1 - \la} - \f{\la}{1 - \la}\\
&=& \f{1 + (e(R) - 2) \la }{1 - \la} .
\eeqn
Now let $\ell(mI^n/xI^n)=0$ for some $n$. Let $s=min\{n|mI^n=xI^n\}.$
Hence by Lemma~ \ref{length} $\ell(mI^n/xI^n)=1$ for all $n=1, 2,
\ldots, s-1.$ Therefore
\beqn
  H(F(I), \la) 
&=& \f{e(R)}{1 - \la} - \ell \left( \f{m}{xR} \right) 
  - \sum_{n = 1}^{s-1} \ell \left( \f{mI^n}{xI^n} \right) \la^n\\
&=& \f{e(R)}{1 - \la} - (e(R) - 1) 
   - \la (1 + \la + \ldots + \la^{s-2})\\
&=& \f{1+ (e(R) - 2) \la + \la^s}{1 - \la} .
\eeqn


\section{The main theorem in dimension two}
The techniques of mixed multiplicities can be illustrated only in rings of
dimension two or higher. The main theorem will be proved by induction on
the dimension of $R.$ Throughout this section $(R,m)$ will
denote a \CM local ring of dimension two with infinite residue field.

Let $I$ and $J$ be  ideals in $R$. Let  $x,y$ be
elements  in $R$. Consider the sequence 
\beq
\label{exact}
0 \, \lrar \, \f{R}{(I:y) \cap (J:x)}
  \, \sta{\psi}{\lrar} \, \f{R}{I} \oplus \f{R}{J}
  \, \sta{\phi}{\lrar} \, \f{(x,y)}{xI+yJ} 
  \, \lrar \, 0
\eeq
where $\phi(a',b') = (ax+by)'$ and $\psi(r') = ((-ry)', (rx)')$. Here
primes denote  the residue classes. 

\bl
\label{exseq}
 If  $x,y$  a
regular sequence in  $R$, then the sequence~(\ref{exact}) is exact. 
\el

\pf Let $(a',b') \in \mbox{Ker}~\phi$. Then $ax+by = xi+yj$ for some $i
\in I$ and $j \in J$. Hence  $x(a-i) = y (j-b)$. Since $x,y$ is a
regular sequence, there exists $t \in R$ such that $a = i -yt$ and $b
= j + xt$. Therefore $(a', b') = ((-yt)', (xt)')=\psi(t^{'})$ which
implies
that
$\mbox{Ker~}\phi = \mbox{Im~} \psi$. 
It is easy to see that  $\psi$ is injective.

\bl
\label{hser2}
 Let $(x,a)$ be a joint
reduction of the set of ideals $(m,I)$. Then  
\beqn
   H(F(I), \la) 
&=& \f{1 + (e_{1}(m|I) -1) \la}{(1 - \la)^{2}}
  - \f{1}{1 - \la} \sum_{n=1}^{\infty}   
    \ell \left( \f{mI^n}{xI^n + amI^{n-1}} \right)\la^n \\
&&+ \f{1}{1 - \la} \sum_{n=1}^{\infty}  
  \ell \left( \f{(mI^{n-1} : x) \cap (I^n : a)}{I^{n-1}} \right)\la^n.
\eeqn
\el
\pf  Put $y=a$, $J=mI^{n-1}$ and $I = I^n$ in 
 (2) to get 
\beqn
\ell \left( \f{(x,a)}{xI^n + a mI^{n-1}} \right)
- \ell \left( \f{R}{I^n} \right) - \ell \left( \f{R}{mI^{n-1}} \right)
+ \ell \left( \f{R}{(I^n :a) \cap (mI^{n-1} : x)} \right) = 0.
\eeqn
Hence
\beqn
  \ell \left( \f{I^n}{mI^n} \right)
- \ell \left( \f{I^{n-1}}{mI^{n-1}} \right)
= e_{1}(m|I)
- \ell \left( \f{mI^n}{xI^n + amI^{n-1}} \right) 
+ \ell \left( \f{(mI^{n-1} : x) \cap (I^n : a)}{I^{n-1}} \right).
\eeqn
Therefore
\beqn
& & (1 - \la) H(F(I), \la) \\
&=& 1 + \sum_{n=1}^{\infty} 
        \left[ \mu(I^{n}) - \mu(I^{n-1}) \right] \la^n \\
&=& 1 + \sum_{n=1}^{\infty} \left[ e_{1}(m|I)
      - \ell \left( \f{mI^n}{xI^n + amI^{n-1}} \right)
      + \ell \left( \f{(mI^{n-1} : x) \cap (I^n : a)}{I^{n-1}} 
                            \right) \right] \la^n\\
&=&     \f{1 + (e_{1}(m|I) -1) \la}{1 - \la} 
      - \sum_{n=1}^{\infty} \left[ 
        \ell \left( \f{mI^n}{xI^n + amI^{n-1}} \right)
      - \ell \left( \f{(mI^{n-1} : x) \cap (I^n : a)}{I^{n-1}} 
                    \right) \right] \la^n.
\eeqn

\bt
\label{maindim2}
   Suppose $\mu(I) =
     e_{1}(m|I)$ and $ \gamma(I) \geq 1$. Then

    \beqn
         H(F(I), \la) 
         = \left\{ 
                  \begin{array}{ll}
        \displaystyle \f{1+ (\mu(I) - 2)\la}{(1 - \la)^2} & 
                       \mbox{ if}\;\; r(m|I)=\infty, \,\, 
                          \\ \\
        \displaystyle \f{1 + \la (\mu(I) - 2) + \la^s}{(1 - \la)^2}& 
        \mbox{if } r(m|I)=s.\\ 
                  \end{array} 
          \right.
     \eeqn
\et

\pf By Lemma \ref{joint}  there exists   a joint
reduction    $(x,a)$  of $(m,I)$ such that $a^o$ is a nonzerodivisor in
$G(I)$.  Hence $(I^n :a) = I^{n-1}$ for
all $n \geq 0$. By
Lemma~\ref{length}, $\ell (mI^n/ xI^n + amI^{n-1}) \leq 1$ for all $n
 \geq 1$. If $\ell (mI^n/ xI^n + amI^{n-1}) = 1$ for all $n
\geq 1$ then by Lemma~\ref{hser2}, 
\beqn
H(F(I), \la)
= \f{1 + (e_{1}(m|I) -1) \la}{(1 - \la)^{2}}
  -  \f{1}{1 - \la} \sum_{n=1}^{\infty} \la^n
= \f{1 + (e_{1}(m|I) -2) \la}{(1 - \la)^{2}}.
\eeqn
If $r(m|I)=s$ then  by Lemma~\ref{hser2}, 
$$
H(F(I), \la)
= \f{1 + (e_{1}(m|I) -1) \la}{(1 - \la)^{2}}
  -  \f{1}{1 - \la} \sum_{n=1}^{s-1} \la^n 
= \f{1 + (e_{1}(m|I) -2) \la + \la^s}{(1 - \la)^{2}}.
$$

 
\section{The main theorem in  dimension $\geq 3$}

In this section we prove the main theorem in \CM local rings of dimension
at least 3. We do this by going modulo a regular element whose initial
form in $F(I)$ and $G(I)$ is simultaneously regular. This preserves
all the hypotheses and we invoke the result in dimension 2. 

\bl
\label{zerodivisor}
Let $(R,m)$ be a local ring. Let $I$ be an ideal of $R$ and $a \in
I\setminus mI.$ 
 Then  $a^{\star}$ is a nonzerodivisor
in $F(I)$ if and only if $(mI^{n+1} : a) \cap I^n = mI^n$ for all $n \geq
0.$
\el

\pf 
Suppose $a^{\star}$ is a nonzerodivisor in $F(I).$ If $b \in (mI^{n+1} :
a) \cap I^n$ then $ba \in mI^{n+1}.$ Thus $b^{\star} a^{\star}=0.$ Since
$a^{\star}$ is a nonzerodivisor, $b^{\star}=0.$ Hence $b \in mI^n.$

Conversely suppose $(mI^{n+1} : a) \cap I^n = mI^n$ for all $n \geq
0.$ As $0:a^{\star}$ is a homogeneous ideal, it is enough to show that
it contains  no  nonzero homogeneous element. Let $b \in I^n$ and 
$b^{\star} a^{\star}=0.$ Then $ba \in mI^{n+1}.$ Hence $b \in (mI^{n+1} :
a) \cap I^n = mI^n,$ which implies that $b^{\star}=0.$

\bl
\label{fiberpreserve}
Let $I$ be an ideal of a local ring $(R,m)$ and $a \in I \sms mI$
be a nonzerodivisor. 
 Suppose $a^{o}$ is a nonzerodivisor in $G(I)$ and $b \in I^n.$ Then  the
map 
$$ 
\phi : F(I)/a^{\star}F(I) \longrightarrow F(I/aR), \hspace*{.25cm} 
\phi((b+mI^n)^{'})=b+(mI^n+aR)
$$
 is an isomorphism.
\el

\pf
As $a^{o}$ is a nonzerodivisor in $G(I)$, $aR \cap I^n=aI^{n-1}$
for all $n \geq 0.$ Therefore $mI^n + aI^{n-1}=mI^n+(aR \cap I^n).$ Thus
$$F(I/aR) \simeq  \bigoplus_{n=0}^{\infty} (I^n+aR)/(mI^n+aR ) 
        \simeq  \bigoplus_{n=0}^{\infty}  I^n/(mI^n+(I^n \cap aR)) 
        \simeq  F(I)/a^{\star}F(I).
$$

\bl
\label{Hilbertseriespreserve}
As per the notation in the above lemma, let $a^{\star}$ and $a^{o}$
be nonzerodivisors in $G(I)$ and $F(I)$ respectively. Then

$$
H(F(I/aR), \lambda) = (1- \lambda )H(F(I), \lambda).
$$
\el

\pf This is clear by the isomorphism in Lemma  \ref{fiberpreserve}.

\bl
\label{indexpreserve}
 Let 
$J=(x, a_1, a_2, \ldots, a_{d-1})$ be a joint reduction of 
 $(m|I^{[d-1]})$. Suppose
$r_J(m|I)$
is finite. If
$a{_1}^{\star}$ is a nonzerodivisor in $F(I)$ and $a_1^0$ is a
nonzerodivisor in $G(I)$ then
$r_J(m|I)=
r_{J'}\left(m'|I'\right)$. Here ' denotes the residue class modulo
$(a_1).$  
\el

\pf 
Suppose $mI^n+ a_1R=xI^n+( a_2, \ldots, a_{d-1})mI^{n-1}+a_1R .$ Then
for any $z \in mI^n$ there exist $ b \in I^n ;\,\, c_2, c_3,
\ldots, c_{d-1} \in mI^{n-1}$ and $ r \in R$ such that
$
z=xb + a_2c_2 + \ldots + a_{d-1}c_{d-1} +a_1r.$

 Hence $ra_1 \in (mI^n \cap a_1R) \subset I^n \cap
a_1R=a_1I^{n-1}$, as
$a_1^0$ is a nonzerodivisor. Thus $r \in I^{n-1}.$ As  $a^{\star}$ is a
nonzerodivisor ,
$r \in (mI^n :a_1) \cap I^{n-1} = mI^{n-1}.$ Thus $ z  \in xI^n+( a_1,
a_2,
\ldots, a_{d-1})mI^{n-1}.$

\bt
\label{maintheorem}
Let $(R,m)$ be a \CM local ring of dimension $d \geq 2.$ Let $I$ be an
$m-$primary ideal of $R$ with almost minimal mixed multiplicity. 
If $\gamma(I) \geq d-1$ and   $\phi(I) \geq d-2$ then

\beqn
H(F(I), \la) 
= \left\{ 
         \begin{array}{ll}
         \displaystyle \f{1+ (\mu(I)-d)\la}{(1 - \la)^d} 
        & \mbox{if} \,\,  r(m|I)=\infty,  \\ \\
          \displaystyle \f{1 + (\mu(I)-d ) \la  + \la^s}{(1 - \la)^d}
          & \mbox{if}\,\, r(m|I)=s. \\ 
\end{array} 
\right.
\eeqn

\et

\pf
 Apply induction on $d.$ The $d=2$ case has been proved in Theorem 
\ref{hser2}. Suppose 
$d \geq 3.$ By Lemma \ref{joint}
there exists a joint reduction $(x, a_1, a_2, \ldots, a_{d-1})$ 
 of $(m|I^{[d-1]})$ where $a_1^o, a_2^o, \ldots, a_{d-1}^o$
is a regular sequence in $G(I)$ and  $a_1^*, a_2^*, \ldots, a_{d-2}^*$ is 
a regular sequence in $F(I).$
 
Put $J =(a_1, a_2, \ldots, a_{d-1})R,$
  $M=(a_1^o, a_2^o, \ldots,a_{d-2}^o)G(I)$,
$K=(a_1^{\star}, a_2^{\star}, \ldots, a_{d-2}^{\star})F(I)$ 
and
$L=(a_1, a_2, \ldots, a_{d-2})R$. Let $'$ denote images in $R/L.$ By
repeated use of Lemma 
\ref{fiberpreserve},
we
get $F(I)/K \simeq F(I').$ Moreover $e_{d-1}(m|I)=e_1(m'|I')$ by
\cite{kv} and 
$$
\mu(I')=\mu(I)-d+2 = e_{d-1}(m|I)=e_1(m'|I').
$$ 

Thus $I'$ has almost minimal mixed multiplicity. As $a_{1}^{\star},
a{_2}^{\star}, \ldots, a_{d-2}^{\star}$ is a regular
sequence in $F(I)$,
$H(F(I),\la)(1-\la)^{d-2}=H(F(I)/K,\la) $. 
The result follows by Theorem \ref{maindim2}.

\bp
Let $(R,m)$ be a \CM local ring of dimension $d \geq 2.$ Let $I$ be an
$m$-primary ideal of $R$. Suppose that $(x, a_1, \ldots, a_{d-1}) $ 
is a joint reduction of $(m|I^{[d-1]})$ satisfying the conditions (i)
$mI^2=xI^2+JmI$, where $J=(a_1, \ldots, a_{d-1})$ and (ii)  
$a_1^{\star}, \ldots, a_{d-1}^{\star}$ is a regular sequence in $G(I).$
Then depth $F(I) \geq d-1.$
\ep
\pf By Theorem 2.8 of \cite{cz1}, it is enough to show that for all
$n \geq 2$, $J\cap mI^n=mJI^{n-1}.$ But $mI^n=xI^n+J^{n-1}mI$ for all $n
\geq 2$. Hence  $ J\cap mI^n=(J\cap xI^n) + J^{n-1}mI.$
It remains to show that $J\cap xI^n \subseteq mJI^{n-1}.$ Let $y\in
J\cap xI^n$. Then $y=xi$ for some $i \in I^n.$ Hence $i \in (J:x)\cap
I^{n}=JI^{n-1}$. Hence $y \in mJI^{n-1}.$


\section{Examples}

\bex
{\em Let $k$ be a field and $R=k[[x,y,z]]$ be the power series
ring. Put $I=m^3$ where $m$ is the unique maximal ideal of $R$. Then
$r(I)=2,\;\; \mu(I)=10$ and $e_2(m|I)=9.$ Thus $\mu(I)=e_2(m|I)+d-2.$
It is easily seen that the fiber cone and the associated graded rings
 of $I$ are \CM. It is easy to see that $r(m|m^3)=2.$ The Hilbert series
of
$F(I)$ is

$$  H(F(I),\la)
 =\sum_{n \geq 0} \binom{ 3n+2 }{ 2} \la^n
 =\sum_{n \geq 0} \left[9\binom{ n+2}{ 2}-9\binom{n+1}{ 1}+1 
\right] \la^n 
 =\f{1+7 \la + \la^2}{(1-\la)^3} 
$$
}           
\eex

\bex
{\em This example shows that for the main theorem in dimension $2$ we need
the depth hypothesis of $G(I).$ Let $R =k [x,y]_m$, where $k$ is an
infinite field, 
$x$ and  $y$ are indeterminates and $m = (x,y)$. Let $I = (x^4, x^3y,
xy^3, y^4)$. Then $\mu(I)=4=e_1(m|I)$ and 
$F(I) =k [x^4, x^3y, xy^3, y^4]$.  Note that for all $n \geq 2$, 
$I^n = m^{4n}$. Hence

$$
H(F(I), \la)
= 1 + 4 \la + \sum_{n=2}^{\infty} (4n+1) \la^n
= \f{1 + 2 \la + 2 \la^2 - \la^3}{(1 - \la)^2}.
$$
As the numerator of the Hilbert series has a negative coefficient, $F(I)$
is not Cohen-Macaulay. Hence depth $F(I)=1.$  Since $G(I)_{+} \subseteq
\mbox{ann}(x^2y^2)^{o}$, depth $G(I)=0.$ }   
\eex

\bex
{\em Let $t$ be an indeterminate and $R=k[[t^4, t^5, t^6, t^7]].$
Consider the
ideal $I=(t^4, t^5, t^6)$. Then $(t^4)I^2=I^3$  and $\mu(I) = e(R)-1.$ As 
$$t^7I \subset I^2=(t^8, t^9, t^{10}, t^{11})= t^4(t^4, t^5, t^6, t^7)$$
the initial form of $t^7$ in $G(I)$
is a zerodivisor. Thus the depth of $G(I)$ is zero. It is easy to see that
$r(m|I)=2.$ The Hilbert series of
$F(I)$ is 
$$H(F(I),\la)= 1+3 \la + \sum_{n\geq 2} 4\la^n +\ldots =
\f{1+2\la+\la^2}{(1-\la)}.$$

\noindent
It follows from \cite{drv}  that $F(I)$ is Cohen-Macaulay. }
\eex

\bex{\em
 In this example we answer the question 3.7 of \cite{dgh} partly.
The question asks whether the fiber cones of all $m$-primary ideals
in a one dimensional  \CM local ring of multiplicity 3 are \CM ?
The answer is no. By Proposition 3.5 of \cite{dgh} it is enough to
consider ideals I minimally generated by 3 elements. By Proposition 
2.3 of \cite{drv}, it follows that for such ideals $F(I)$ is \CM if and
only if $r(I)\leq 1.$   Consider the semigroup ring
 $R=k[[t^3, t^7, t^{11}]]$ and the ideal $I=(t^6, t^7, t^{11}).$ Then
$e(R)=\mu(I)=3.$ However $r(I)=2.$ Thus $F(I)$ is not \CM.
}   
\eex

\end{document}